\documentclass[12pt]{article}
\usepackage{amssymb}
\usepackage{epsfig}
\usepackage{color}
\usepackage{amsfonts,amsthm,amsmath,amssymb,amscd,graphicx}

\begin{document}
\begin{center}
{\large\bf On three types of dynamics, and the notion of attractor.}

\vspace{12pt}

{\bf S.V.Gonchenko}$^1$ and {\bf D.Turaev}$^{1,2}$\label{Author}
\vspace{6pt}

{\tiny
$^1$ Lobachevsky University of Nizhny Novgorod, Russia; E-mail: gonchenko@pochta.ru

$^2$ Mathematics Department, Imperial College, London; E-mail: dturaev@imperial.ac.uk}

\begin{abstract}
We propose a theoretical framework for an explanation of the numerically discovered phenomenon of the attractor-repeller merger. We identify regimes which are observed in dynamical systems with attractors as defined in a work by Ruelle and show that these attractors can be of three different types. The first two types correspond to the well-known types of chaotic behavior - conservative and dissipative, while the attractors of the third type, the reversible cores, provide a new type of chaos, the so-called mixed dynamics, characterized by the inseparability of dissipative and conservative regimes. We prove that every elliptic orbit of a generic non-conservative time-reversible system is a reversible core. We also prove that
a generic reversible system with an elliptic orbit is universal, i.e., it displays dynamics of maximum possible richness and complexity.
\end{abstract}

\end{center}

\section{Introduction}{\em Attractor-repeller merger.}
When we speak of dynamical chaos, we usually mean one of two quite different types of dynamics. In Hamiltonian systems, we have conservative chaos -- something like a ``chaotic sea'' with elliptic islands inside. Chaos in dissipative systems is quite different and is associated with strange attractors. Our goal in this paper is to attract attention to one more type of chaos, the third one, which was called ``mixed dynamics'' in \cite{GStS02,DGGLS13}. This type of behavior is characterized
by {\em inseparability} of attractors, repellers, and conservative elements in the phase space \cite{GST97}.

\begin{figure}[ht]
\includegraphics[width=12cm]{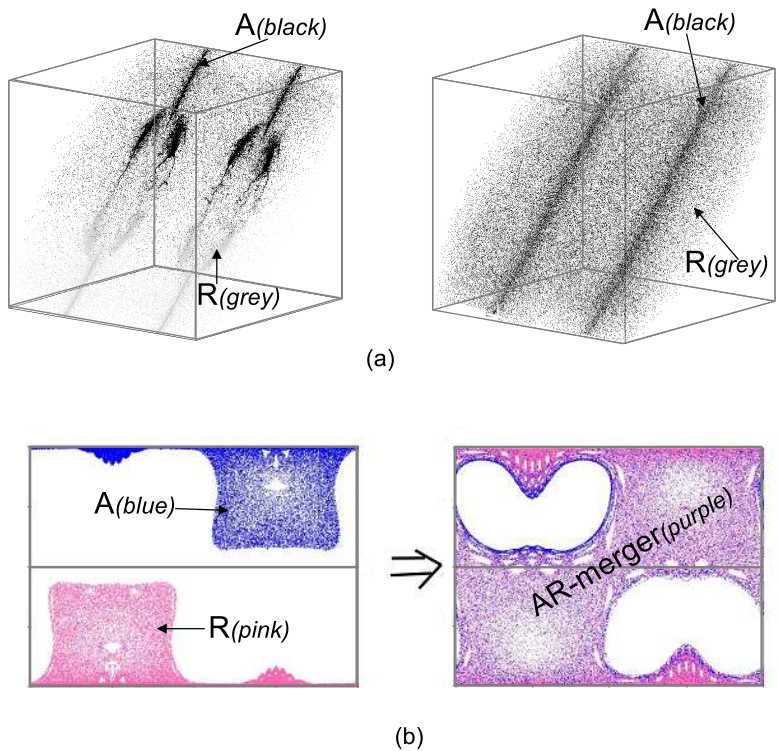}
\caption{{\footnotesize Examples of attractor-repeller merger for the Poincar\'e map of (a),(b) a model of the Celtic stone, \cite{GGK13}, and (c),(d) Chaplygin top (rubber body), \cite{K13}. Here, the numerically obtained attractor ($A$) and repeller ($R$) are shown for different values of the energy of the system.}}
\label{recsc}
\end{figure}

In order to have both attractors and repellers, the system must contract the phase volume somewhere and somewhere expand it. For example, any diffeomorphism of a compact phase space will have an attractor and a repeller, unless the whole phase space is a chain-transitive set (see
the definition of the chain transitivity e.g. in \cite{AnBr85} and Section 2). Attractors and repellers may be separated from each other, like in Morse-Smale systems. However, it was recently established for many examples \cite{PikTop02,GGK13,K13,K16,GGKT17},
that when parameters of a system are varied, the numerically obtained attractor and repeller may {\em collide} and
start to occupy approximately the same part of the phase space. Moreover, a further change of parameters does
not seem to break the attractor-repeller merger. An example of such behavior in models of a Celtic stone and a balanced
rubber ball with that roll on a horizontal plane is shown in Fig.\ref{recsc}.

A theorem by Conley \cite{Co78} establishes the existence of a Lyapunov function, non-increasing along the orbits of the
system, which takes its maxima on repellers and minima on attractors. The attractor-repeller collision means that this
function must have very degenerate critical points, so the intuition based on the Conley theorem might suggest that this
phenomenon should be very exotic. However, in reality, the attractor-repeller merger appears quite robustly
and persists for significant regions of parameter values in the models where it was detected \cite{GGK13,GGKT17}.
As we will argue below, it is a generic phenomenon for the non-conservative time-reversible systems and, more generally,
for systems belonging to the so-called absolute Newhouse domain \cite{T10,T15}.

A systematic study of this phenomenon requires a precise definition of what one means by attractor or repeller.
In this paper we discuss two closely related definitions, going back to the Ruelle's work \cite{Ru81}, which we think
are most relevant in describing numerically obtained portraits of attractors/repellers.
Accordingly, we propose two paradigmatic models (see subsections ``reversible core'' and ``full attractor'' below) for
the attractor-repeller collision, the second less restrictive than the first one,
and also discuss basic dynamical phenomena associated with them in the context of time-reversible systems.
The models are different but not mutually exclusive. In fact, the dynamics underlying the attractor-repeller collision
is so extremely rich that both our models appear to be applicable at the same time.

\bigskip

{\em Reversible core.}
The first model employs the notion of attractor as defined in the works of Ruelle \cite{Ru81} and Hurley \cite{Hu}. Both authors attribute important ideas to Conley \cite{Co78}, so we use the term Conley-Ruelle-Hurley (CRH) attractor -- a stable (with respect to permanently acting perturbations) and chain-transitive closed invariant set. A CRH-repeller is a CRH-attractor in the reverse time; see precise definitions in Section \ref{arc}. We
also show that only three types of CRH-attractors
are possible for a homeomorphism $f$ of a compact, connected, separable metric space
$\mathcal M$:\\
$\bullet$ {\em conservative}
-- when the whole phase space $\mathcal M$ is a chain-transitive set and, hence, the unique CRH-attractor and repeller (this includes e.g. the case of volume-preserving maps on compact manifolds);\\
$\bullet$ {\em dissipative}
-- when an $\varepsilon$-orbit from outside converges to a neighbourhood of the CRH-attractor,\\
$\bullet$ {\em mixed} or
a {\em reversible core} --  a CRH-attractor which is, at the same time, a CRH-repeller, so it retains both forward and backward orbits of $f$ in its small neighborhood (see Sec.~\ref{sec:rc} and Th.~2).

The reversible core should be distinguished from a dissipative attractor, as the reversible core does not attract any orbit. The dynamics here is not exactly conservative either, since the reversible core, as we show in Section \ref{arc} (Theorem 1), is always a limit of a sequence of attractors (and the limit of a sequence of repellers as well). Therefore, we can associate this type of stable sets with the third, mixed type of dynamical behavior.

\begin{figure}[ht]
\includegraphics[width=14cm]{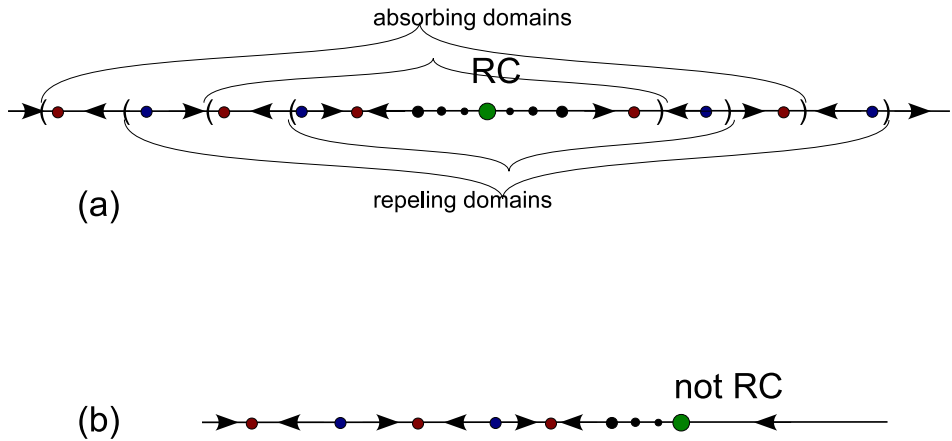}
\caption{{\footnotesize Simple examples of a CRH-attractor that is (a) a reversible core and (b) not reversible core, though it is a limit of CRH-repellers.}}
\label{revcor}
\end{figure}

We stress that there are no further possibilities in this scheme (no ``fourth'' type of chaos). The CRH attractors and repellers
cannot have non-trivial intersections: a CRH-attractor either does not intersect with any CRH-repeller, or it coincides with
some. In the latter case such CRH-attractor is a reversible core. It is easy to construct trivial examples of a reversible core (see Fig.\ref{revcor}a)\footnote{It is also easy to construct examples where a CRH-attractor is a limit of CRH-repellers but it is not a reversible
core (see Fig.\ref{revcor}b). Non-trivial generic examples of such situation can be found in \cite{Bo13}.}. We,
however, also provide a non-trivial example in Section \ref{el2} (see Theorem 3) and, in fact, show that the reversible cores are present {\em generically} in non-conservative time-reversible systems and, hence, are relevant beyond abstract schemes of topological dynamics. Numerical simulations with several models of mechanics provide a direct evidence
for the possible existence of a reversible core in these models: the attractor and repeller in the left part of Fig.~\ref{recsc} are separated, while in the right part they  appear to really coincide!

In theory, one could perform a quite straightforward procedure for the numerical detection of a reversible core: if a numerically obtained attractor does not separate from a numerically obtained repeller with the increase of the accuracy of computation, then this is a reversible core. As we mentioned, a reversible core is always a limit of
an infinite sequence of attractors  (see Theorem 1), so this observation gives a numerical criterion for the coexistence of infinitely many attractors in the phase space. Other known criteria for such phenomenon are based on completely different ideas \cite{Ne74,GST97,LS04,GST09,T10,GO13,DGGLS13,GLRT14,T15}.

\bigskip

{\em Full attractor.}
In reality, however, the computations are rarely repeated many times with ever increasing accuracy.
It is difficult, and may be unnecessary, to distinguish whether the attractor and repeller coincide exactly or whether they
are just very close to each other: in both cases one has a right to speak about a mix of attracting and repelling dynamics.

This leads to a more relaxed idea of the intersection of attractor and repeller, based on a different notion of attractor. We define {\em the full attractor} of the map $f$ as the closure of the union of all its CRH-attractors and {\em the full Ruelle attractor}\footnote{Ruelle did not call this set attractor, but he introduced it for the study of the behavior of epsilon-orbits -- this is the set $A^*$ in Corollary 5.6 of \cite{Ru81}.}
as the prolongation of this set (the prolongation of a set $A$ is the set of all points attainable from $A$ by $\varepsilon$-orbits for arbitrarily small $\varepsilon$; the precise definitions of the full attractor is given in Section \ref{arc}).

While closely related to the behavior of the system subject to a small bounded noise, i.e, to the pictures of the dynamics obtained from numerical or other experiments, such defined attractors are purely dynamical objects. Namely, these are closed invariant sets of the map $f$ and they are preserved by the topological conjugacy: if maps
$f$ and $g$ are conjugate by a homeomorphism $h$, then the full attractor of $f$ is taken by $h$ to the full attractor of $g$, and the full Ruelle attractor of $f$ is taken to the full Ruelle attractor of $g$. The full repeller and the full Ruelle repeller
of $f$ are defined as the full and, respectively, full Ruelle attractors of $f^{-1}$.

If the entire phase space is chain-transitive, then the attractor and repeller are equal to the whole of the phase space, as we already mentioned; we consider such dynamics conservative from the topological point of view. When the full Ruelle attractor and the full Ruelle repeller of $f$ do not intersect, we say that the global dynamics of $f$ is dissipative. In the last remaining case, where the phase space is not chain-transitive but the full Ruelle attractor
and repeller have a non-empty intersection, we say that the global dynamics of $f$ is mixed. As any reversible core is simultaneously a CRH-attractor and a CRH-repeller, it belongs to both the full attractor and the full repeller, so the existence of reversible cores implies the mixed dynamics of $f$. However, there are
more general possibilities for the mixed dynamics, as discussed below.

\bigskip

\begin{figure}[ht]
\centerline{\includegraphics[width=8cm]{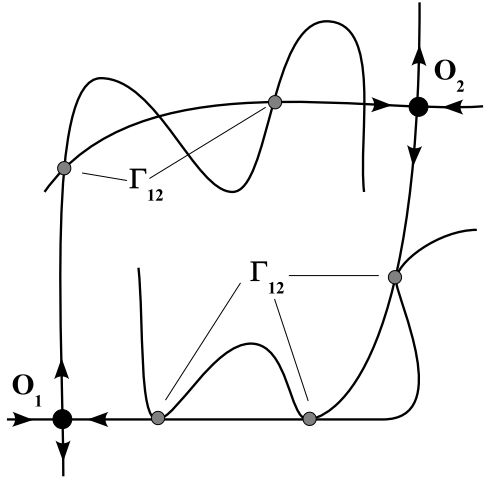}}
\caption{
{\footnotesize
A non-transverse heteroclinic cycle for a two-dimensional diffeomorphism: two saddle fixed points, $O_1$ and $O_2$,  and two heteroclinic orbits, $\Gamma_{12}$ and $\Gamma_{21}$, such that $W^u(O_1)$ and $W^s(O_2)$ intersect transversely at the points of  $\Gamma_{12}$  and $W^u(O_2)$ and $W^s(O_1)$ have a quadratic tangency at the points of  $\Gamma_{21}$. When
$(|J(O_1)|-1)(|J(O_2)|-1)<0$, bifurcations of this cycle lead the system
into the absolute Newhouse domain and create a non-destructible intersection of the full attractor and full repeller.}}
\label{cont}
\end{figure}

{\em Heteroclinic cycles at the intersection of attractor and repeller.}
In the retrospect, the phenomenon of a non-removable intersection of the full attractor and repeller was discovered in \cite{GST97}, where we proved that the closure of the set of asymptotically stable periodic orbits (sinks) and the closure of
the set of repelling periodic orbits (sources) may have a persistently non-empty intersection. Namely, for a class of
heteroclinic cycles shown in Fig.~\ref{cont}, in any continuous family of smooth two-dimensional maps for which the heteroclinic tangency splits there exist open regions (Newhouse regions) in the parameter
space where a residual set of parameter values corresponds to a non-empty intersection of the closure of sinks and the
closure of sources.

The main property of the heteroclinic cycles that produce mixed dynamics is that they contain two
saddles such that the map is area-contracting near one of the saddles and area-expanding near the other (the maps having
cycles with this property are dense in the absolute Newhouse domain of \cite{T10,T15}).

\begin{figure}[h!]
\centering
\includegraphics[width=14cm]{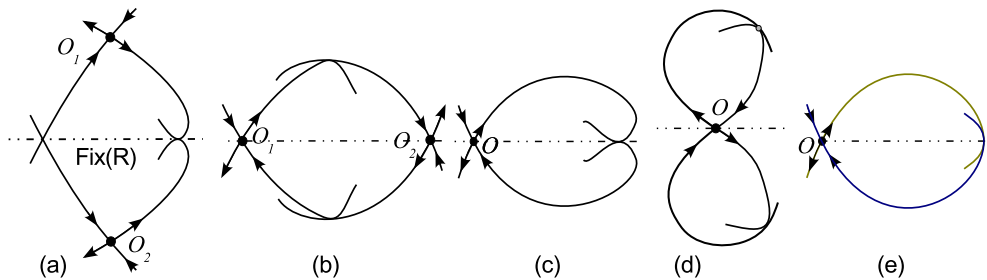}
\caption{{\footnotesize Examples of codimension-1 bifurcations of homoclinic and heteroclinic tangencies in two-dimensional reversible maps. Maps with nontransverse heteroclinic cycles are shown in (a) and (b): here (a) $O_1 = g(O_1)$ and $J(O_1)=J(O_2)^{-1}<1$, (b) $J(O_1)=J(O_2) =1$. Maps with homoclinic tangencies are shown in (c)--(e): here the point $O$ is symmetric in all cases; the homoclinic orbit is symmetric in cases (c) (quadratic tangency) and (e) (cubic tangency); (d) gives an example of reversible map with a symmetric pair of quadratic homoclinic tangencies to $O$.}}\label{5typs}
\end{figure}

Heteroclinic cycles with this property naturally appear in reversible diffeomorphisms (a diffeomorphism $f$ is reversible if it
is conjugate to its own inverse $f^{-1}$ by means of a certain involution $g$). An example of a heteroclinic cycle of the above described type for a reversible diffeomorphism is shown in Fig.~\ref{5typs}a. Note a symmetry in this picture: the involution $g: ( x\to x, y\to -y)$ takes the saddle fixed point $O_1$ to the saddle fixed point $O_2$, the invariant
manifold $W^u(O_1)$ to $W^s(O_2)$, and $W^u(O_2)$ to $W^s(O_1)$. In general, the Jacobian $J$ of the
map at $O_1$ can be arbitrary; for instance, let $|J(O_1)|<1$. Since $f^{-1}$
near $O_2$ is smoothly conjugate to $f$ near $O_1$, it follows that $|J(O_2)|=|J^{-1}(O_1)|>1$. It was shown in \cite{LS04} that if the map $f$ is embedded into a family of reversible maps for which the non-transverse heteroclinic
in such cycle splits, then for generic values of parameters from the corresponding Newhouse regions there exist infinitely many periodic sinks, sources, saddles, and {\em elliptic periodic orbits}. Moreover, the closure of the set of each of these types of periodic orbits contains the points $O_1$ and $O_2$, i.e. {\em attractors, repellers, saddles, and elliptic orbits are
inseparable from each other}.

Similar results for reversible maps with another type of heteroclinic
cycle (see Fig.~\ref{5typs}b), are obtained in \cite{DGGLS13}. Unlike the previous case, the saddle points $O_1$
and $O_2$ belong to the line of the fixed points of the involution $g$, hence $J(O_1)=J(O_2) =1$. The involution $g$
takes $W^u(O_1)$ to $W^s(O_1)$, and $W^u(O_2)$ to $W^s(O_2)$, and the map $f$ has a symmetric pair of non-transverse heteroclinic orbits.
In contrast to the case of Fig.~\ref{5typs}a, the map $f$ near the fixed points $O_{1,2}$ is conservative. However,
the conservativity is violated in this situation near the heteroclinic tangencies -- the maps along heteroclinic tangencies have, in general,
a non-constant Jacobian\footnote{Corresponding explicit conditions of general position
are given in \cite{DGGLS13}.}. As shown in \cite{DGGLS13}, bifurcations of such heteroclinic cycle also
lead to the the reversible mixed dynamics (i.e., an unbreakable intersection of the closures of the sets of attractors, repellers, saddles, and eliptic orbits). The same phenomenon takes place at the bifurcations of a symmetric pair of homoclinic tangencies as in Fig.~\ref{5typs}d (here $g(O)=O$ and $g(W^u(O)) = W^s(O)$), see \cite{DGGL15}.
We also plan to prove the existence of the absolute Newhouse intervals
in one-parameter families of reversible maps, which unfold symmetric quadratic and cubic homoclinic tangencies, as in Figs.~\ref{5typs}c and~\ref{5typs}e, respectively.

\bigskip

{\em Richness of reversible mixed dynamics.}
In this paper we undertake a further investigation of the mixed dynamics in reversible maps. We show (Section \ref{el2}) that every
elliptic periodic orbit of a generic two-dimensional reversible diffeomorphism is a reversible core. We also show (Section \ref{eln}) that
every elliptic periodic orbit of a generic $n$-dimensional ($n\geq 2$) reversible diffeomorphism is {\em a limit
of a sequence of uniformly-hyperbolic attractors and repellers of all topological types possible in an $n$-dimensional ball}.

We recall basic definitions from the theory of non-conservative reversible systems in Section \ref{eln}. Such systems
are known to appear in various applications. In particular, non-conservative time-reversible dynamical systems
are natural models for mechanical systems with nonholonomic constraints (e.g. the examples of dynamics shown in Fig. \ref{recsc} are obtained for
reversible systems of this type). The existence of elliptic periodic orbits is a characteristic property of reversible maps,
for which the dimension of the
set $\mbox{Fix}(g)$ of the fixed points of the involution $g$ is half of the dimension of the phase space, or higher. In particular, they emerge
in various homoclinic bifurcations \cite{LS04,DGGLS13,DGGL15,FT96}. According to the Reversible Mixed Dynamics conjecture of \cite{DGGLS13}, the mixed
dynamics emerging at the most typical homoclinic bifurcations of reversible systems must always include a large number of elliptic periodic orbits.

It is known \cite{Sevr} that the dynamics of a reversible map near an elliptic periodic orbit is, to a large extent, conservative: a significant portion
of the phase space in a neighborhood of such orbit is filled by KAM-tori. However, the dynamics in the resonant zones between the
KAM-tori does not, in general, need to be conservative \cite{GLRT14}. Results of Section \ref{eln} show that the dynamics near
a typical elliptic orbit is {\em universal} in the sense of \cite{T03,GST06,GST07,T10,T15}, i.e., the iterations of the Poincar\'e map
in the resonant zones near the elliptic orbit approximate, with arbitrarily good precision, every dynamics possible in the given dimension of the
phase space. Thus, on one hand non-conservative reversible systems are of immediate importance for applications. On the other hand their
dynamics is an example of an ultimate richness and complexity: any generic non-conservative reversible system with an elliptic point exhibits {\em all robust dynamical phenomena possible}.

\section{Attractors, repellers, and a reversible core}\label{arc}

In several papers \cite{GStS02,GSSt06,G10,G13,DGGLS13,GLRT14,G16,GGKT17} a new, third type of chaotic dynamics was identified -- the so-called ``mixed dynamics'' characterized by the attractor-repeller merger. Below we propose a scheme which could formalize this idea.

\subsection{Definitions of the attractor.}
We start with reminding definitions and simple facts from topological dynamics. Consider a homeomorphism
$f$ of a compact, separable, metric space $\mathcal M$.
A sequence of points $x_1,\dots,x_N$ is called {\em an $\varepsilon$-orbit} of the map $f$ if $dist(f(x_j), x_{j+1})<\varepsilon$ for all $j=1, \dots, N-1$.
We will say that the $\varepsilon$-orbit $x_1,\dots,x_N$ {\em connects} the point $x_1$ to $x_N$ and that $x_N$ is {\em attainable} from $x_1$ by an $\varepsilon$-orbit of the length $N$.

A closed invariant set $\Lambda$ is called {\em chain-transitive} if for every $\varepsilon>0$ and every two points $x\in\Lambda$ and $y\in\Lambda$ there exists
an $\varepsilon$-orbit that lies in $\Lambda$ and connects $x$ and $y$.
If
two chain-transitive sets $\Lambda_1$ and $\Lambda_2$ have a non-empty intersection, their union $\Lambda_1\cup\Lambda_2$ is also a chain-transitive set.

An open set $U\subset \mathcal M$ is {\em an absorbing domain} if $f(cl(U))\subset U$.  An important fact is that
\begin{itemize}
\item
the set $D_{K,\varepsilon}(x)$ of all points attainable from $x$ by $\varepsilon$-orbits of the length $N\geq K$ is {\em always an absorbing domain}.
\end{itemize}

Indeed, this set is obviously open, so we just need to show
that for any point $z$ in the closure of $D_{K,\varepsilon}(x)$ there is an $\varepsilon$-orbit of length at least $K$ that connects $x$ and $f(z)$.
By definition, $z$ can be approximated arbitrarily well by an end point $x_N$ of some $\varepsilon$-orbit $x_1=x,\dots,x_N$, where $N\geq K$. Let $dist(x_N,z)<\delta$
where $\delta$ is such that the images of any two $\delta$-close points by $f$ lie at a distance smaller than $\varepsilon$ from
each other. Then $dist(f(x_N),f(z))<\varepsilon$, i.e. $x_1,\dots,x_N,f(z)$ is an $\varepsilon$-orbit that connects $x_1=x$ to $f(z)$, i.e. $f(z)\in D_{K,\varepsilon}(x)$.

The only absorbing domain that a {\em volume-preserving} map of a compact manifold $\mathcal M$ can have is $\mathcal M$ itself. Therefore, for any volume-preserving map
$D_{K,\varepsilon}(x)=\mathcal M$ for every $x$, i.e., every two points of $\mathcal M$ are connected by an $\varepsilon$-orbit for any $\varepsilon>0$. Thus,
$\mathcal M$ is a chain-transitive set in this case. There are other, not volume-preserving examples. For instance, for any homeomorphism
which is topologically conjugate to a volume-preserving map (like e.g. any Anosov map on a torus) the whole phase space is chain-transitive (if the phase space is compact);
another example is given by the map $\displaystyle\varphi\mapsto \varphi+\sin^2\frac{\varphi}{2}$ of a circle, which has a single semi-stable fixed point, etc.

We will not consider  anymore the case where the whole phase space is chain-transitive. Thus, we assume that some points of $\mathcal M$ are not {\em chain-recurrent}
(a point $x$ is chain-recurrent if its $\varepsilon$-orbits return to it infinitely many times, i.e. if $x\in D_{K,\varepsilon}(x)$ for all $\varepsilon>0$ and $K>0$).

Denote $D(x)=\bigcap_{K,\varepsilon} D_{K,\varepsilon} (x)$.
This is a closed invariant set and it is an intersection of a family of nested absorbing domains. Such sets are  {\em stable} meaning that given
any $\delta>0$ there exists a sufficiently small $\varepsilon>0$ such that no $\varepsilon$-orbit starting at this set can leave its $\delta$-neighbourhood.
If a stable set $A$ is chain-transitive, it has no proper stable subsets. If not, it contains a point $y$ which is not chain-recurrent, and the set $D(y)$ is a proper
and stable subset of $A$. By (transfinite) induction, one can show that {\em every stable set has a chain-transitive stable subset}, cf. \cite{Ru81,Hu}. Following Ruelle and Hurley,
we will call a {\em chain-transitive and stable set} an attractor of the map $f$, or {\em Conley-Ruelle-Hurley-attractor} (CRH-attractor).

We identify observable dynamical regimes with trajectories which stay in a neighborhood of CRH-attractors.
The logic is:  Whenever a certain dynamical process is observed, there is never a guarantee that
the dynamical system, which generates it, is known precisely (for example, when we compute orbits of a given map numerically, the resulting sequence of points is, in fact, an orbit of a slightly different map, due to the round-off). One therefore may claim that the observed regimes are $\varepsilon$-orbits with
a sufficiently small $\varepsilon$.

As a simplified model, one can consider the $\varepsilon$-orbits $\{x_j\}$ of the map $f$ as realizations of a random process such that the deviations $x_{j+1}$ from $f(x_j)$ are {\em independent random variables} $\xi_j$ with probability densities supported in the ball
$\|\xi_j\|\leq \varepsilon$, which are bounded away from zero and continuous in this ball. Then, for a fixed $\varepsilon>0$ it is natural to define the
$\varepsilon$-{\em attractor} of a point $x_0$ as the set $A_\varepsilon(x_0)$ of all points which are $\omega$-limit for $\varepsilon$-orbits of $x_0$ with positive probability. The result does not depend on the choice of the
probability density for $\xi_j$ -- it is a union of finitely many closed sets such that no forward $\varepsilon$-orbit starting in any of these sets can leave it -- so the interior of this set is an absorbing domain, and in each of these absorbing domains every two points are connected by an $\varepsilon$-orbit. It immediately follows that if some point $x^*$ belongs to the intersection of some $\varepsilon_j$-attractors for a sequence $\varepsilon_j\to 0$ as $j\to+\infty$, then $x^*$ belongs to the intersection of a sequence of nested absorbing domains, and every two points in this intersection are connected by $\varepsilon_j$-orbits for all $j$, i.e. $x^*$ belongs to a stable and chain-transitive set -- a CRH-attractor. Thus, every point which is not
in a CRH-attractor stays outside of $\varepsilon$-attractors for all sufficiently small $\varepsilon$.

Another way
to express the same idea is to note that the intersection
$A_0(x)=\bigcap_{\rho>0} \bigcup_{0<\varepsilon\leq \rho} A_{\varepsilon}(x)$ consists only of CRH-attractors (moreover, each of these CRH-attractors is attainable from $x$ by $\varepsilon$-orbits for every $\varepsilon$ -- we will further simply say that it is attainable from $x$). We note that the sets
$\cup_{0<\varepsilon\leq \rho} A_{\varepsilon}(x)$ are not closed, so it is difficult to investigate
the structure of their intersection $A_0(x)$. However, one can show that $A_0(x)$ is non-empty. Moreover,
its closure $\overline{A_0(x)}$ contains {\em all} CRH-attractors attainable from $x$, so it is the closure of the union of all
CRH-attractors, attainable from $x$.

Thus, the following definition makes sense:
\begin{itemize}
\item
{\em an attractor of a point} $x$ is any CRH-attractor attainable from $x$.
\end{itemize}
If the number of such attractors is finite, then their union is the {\em full attractor} of $x$.
In the case of an infinite number of such attractors, there is less certainty in the definition of the full attractor.
One candidate would be {\em the closure of the union of all attractors of} $x$, i.e. the above defined set
$\overline{A_0(x)}$. It is a closed invariant set, but it may be not stable. It is easy to show that the minimal closed stable set that contains $\overline{A_0(x)}$ is the {\em prolongation} of $\overline{A_0(x)}$, i.e. the set of all points which are attainable from $\overline{A_0(x)}$ for every arbitrarily small $\varepsilon$. We call it the {\em full Ruelle attractor} of $x$. Similarly, one defines:
\begin{itemize}
\item
 the {\em full attractor of the map} $f$ is {\em the closure of the union of all its CRH-attractors}, and the {\em full Ruelle attractor of} $f$ is {\em the prolongation of the full attractor}. In the same way one produces definitions of repellers as attractors for the inverse map $f^{-1}$.
\end{itemize}

\subsection{Absolute Newhouse domain.}\label{absnew}
With these definitions, our goal is to investigate how attractors and repellers can intersect. The possibility of intersection of the full attractor and the full repeller was discovered in \cite{GST97}. The main idea of that paper is that a generalization of the classical Newhouse phenomenon to maps which
are not area-contracting
leads to a persistent intersection of the full attractor and full repeller (as we defined them here).

It was shown by Newhouse in \cite{Ne79, Ne74} that a two-dimensional diffeomorphism can have a {\em wild hyperbolic set} -- a zero-dimensional compact transitive hyperbolic set $\Lambda$ whose stable and unstable sets $W^s(\Lambda)$ and $W^u(\Lambda)$ have a tangency which is not removable by any $C^2$-small perturbation. By the definition, maps with wild hyperbolic sets form a $C^2$-open region in the space of two-dimensional diffeomorphisms, we call this region {\em the Newhouse domain}. If we take a map from the Newhouse domain, then every map from its small neighborhood $\mathcal D$ in the space of $C^2$-diffeomorphisms will have a hyperbolic set $\Lambda$ such that $W^s(\Lambda)$ and $W^u(\Lambda)$ are tangent. Since saddle periodic orbits are dense in $\Lambda$ and their stable and unstable manifolds are dense in $W^s(\Lambda)$ and $W^u(\Lambda)$, respectively,  the maps which have a homoclinic tangency
of the stable and unstable manifold of some periodic orbit in $\Lambda$ are $C^r$-dense in $\mathcal D$
(for every $r\geq 2$).

As it was shown first by Gavrilov and Shilnikov in \cite{GaS73}, bifurcations of a homoclinic tangency to a saddle periodic orbit of a two-dimensional map lead to the birth of {\em stable periodic orbits if the Jacobian $J$ of the first-return map at the saddle periodic point is smaller than 1 in the absolute value}. Using this fact (which he established independently), Newhouse
proved in \cite{Ne74} that maps with {\em infinitely many periodic attractors} {\em (stable periodic orbits)}
are {\em dense} ({\em and form a residual subset}) in the part of the Newhouse domain of the space of $C^r$-diffeomorphisms that corresponds to maps which are {\em area-contracting} in a neighborhood of
the wild hyperbolic set. Moreover, for a $C^r$-generic map of this class {\em the closure of the set of stable
periodic orbits contains a wild hyperbolic set} $\Lambda$ {\em and all points homoclinic to} $\Lambda$.

In \cite{GST97} we considered heteroclinic cycles with two saddle periodic orbits, $O_1$ and $O_2$, such that
$|J(O_1)|<1$ and $|J(O_2)|>1$ (where $J(O)$ is the product of the multipliers of $O$, i.e., the determinant of the
derivative of the first-return map at $O$). The cycle contains orbits of the heteroclinic intersection between
$W^u(O_1)$ with $W^s(O_2)$ and $W^u(O_2)$ with $W^s(O_1)$, and we assume that one of these orbits corresponds
to a tangency of the corresponding stable and unstable manifolds. We showed in \cite{GST97} that
a generic splitting of this tangency leads to the creation of a wild hyperbolic set $\Lambda$ which contains
both orbits $O_1$ and $O_2$, i.e., $\Lambda$ contains a pair of saddle periodic orbits such that the map is
area-contracting near one of the orbits and area-expanding near the other. Such wild hyperbolic sets were called in \cite{T15} {\em ultimately wild}; the part of the Newhouse domain which contains maps with ultimately wild sets is called the {\em absolute Newhouse domain} \cite{T10,T15}. Applying Gavrilov-Shilnikov results \cite{GaS73}, we showed in \cite{GST97} that
\begin{itemize}
\item
{\em a generic map in the absolute Newhouse domain has infinitely many coexisting periodic attractors and repellers, and the closure of the set of these attractors and the closure of the set of these repellers both contain the ultimately wild hyperbolic set} $\Lambda$ {\em and all points homoclinic to it}\footnote{For three-dimensional diffeomorphisms, analogous results were established in \cite{GST09,GO10,GO13}; the existence of absolute Newhouse domain in any dimension was shown in \cite{T96}}.
\end{itemize}
For such maps the full attractor and full repeller both contain the wild set $\Lambda$ and its homoclinic points, so they have a non-empty intersection. In spite the intersection, the full Ruelle attractor and the full Ruelle repeller can differ significantly, as the full Ruelle attractor also contain the whole of unstable manifold of $\Lambda$, while the full Ruelle repeller contains the whole of the stable manifold of $\Lambda$.\footnote{Note that in the case of {\em area-preserving} maps on compact surfaces, the closures of stable and unstable manifolds of $\Lambda$ coincide \cite{Tresh1,Tresh2}, and this is in a complete agreement with the fact that the attractor coincides with the repeller in this case.}

We stress that this particular instance of the intersection of the attractor and repeller is, probably, the most basic model for the mixed dynamics in two-dimensional diffeomorphisms. Indeed, if we see something which looks like a chaotic attractor, then it is natural to expect that it contains a hyperbolic set and areas are contracted near this set. Similarly, it is natural to expect that a chaotic repeller contains a hyperbolic set near which areas are expanded by the map. Collision of such sets would involve creation of heteroclinic connections between the saddle periodic orbits belonging to these two sets. This means the existence of heteroclinic tangencies at a sequence of parameter values during the attractor-repeller collision process. Thus, the above described heteroclinic cycles appear, so the maps undergoing the attractor-repeller merger can be with a good degree of certainty be considered as belonging to the absolute Newhouse domain.

Note that the generic ultimately-wild hyperbolic sets serve as the limits of objects much more complicated  than just stable and unstable periodic orbits, It was shown in \cite{GStS02,GSSt06} that they, generically, are the limits of
sequences of stable and unstable closed invariant curves (i.e. quasiperiodic attractors and repellers). In \cite{T15} it was shown that a generic ultimately-wild hyperbolic set is a limit of a sequence of hyperbolic attractors and hyperbolic repellers of all topological types possible for a diffeomorphism of a two-dimensional disc. In fact, it was shown in \cite{T15} that the dynamics of a $C^\infty$-generic two-dimensional diffeomorphism in an arbitrarily small neighborhood of the closure of the set of points homoclinic to
an ultimately wild set is {\em universal}, i.e. iterations of such map provide arbitrarily good $C^r$
approximations to {\em all} orientation-preserving $C^r$-diffeomorphisms of a two-dimensional disc into $R^2$, for every $r$. In this sense, the mixed dynamics near the intersection of the full attractor and a full repeller of a generic diffeomorphism belonging to the absolute Newhouse domain is of maximal possible richness.

\subsection{Reversible core.} \label{sec:rc}
Next, we report a different, previously unnoticed mechanism for the coexistence of infinitely many attractors and repellers and their intersection. It is based on the observation that there are, in fact, two types of the CRH-attractors: those which actually attract something, and those which attract nothing.

Namely,
\begin{itemize}
\item
we call a CRH-attractor $A$ {\em a dissipative attractor} if there is a point $x\not\in A$ such that for any $\varepsilon>0$ there is an $\varepsilon$ orbit that connects $x$ to a point in $A$. Otherwise we will call the CRH-attractor $A$ {\em a reversible core}.
\end{itemize}

It is obvious that every absorbing domain (different from the whole phase space $\mathcal M$) must contain at least one dissipative attractor. However, the reversible cores may also exist (and they can exist $C^\infty$-generically as we show in the next Section). As it is explained next, the existence of the reversible core immediately implies mixed dynamics.

By the definition, given any $\delta>0$ there exists $\varepsilon>0$ such that no $\varepsilon$-orbit that starts at a distance larger than $\delta$
from the reversible core $C$ can end at a point of $C$. Note that the set $U_\varepsilon(C)$ of all points $x$ such that some $\varepsilon$-orbit of $x$
ends in $C$ is open, and it is easy to see that it satisfies $f^{-1}(cl(U))\in U$, i.e. it is an absorbing domain for the inverse map $f^{-1}$. Thus, any neighborhood
of $C$ contains an absorbing domain for the map $f^{-1}$, hence the reversible core is a CRH-attractor for the inverse map as well, i.e., it is also a CRH-repeller. Thus, the reversible core is an intersection of a sequence
of embedded absorbing domains for the map $f$ and absorbing domains for the map $f^{-1}$, see Fig.~\ref{revcor}(b).
Any absorbing domain that encloses $C$ must contain at least one dissipative attractor, and any absorbing domain for the inverse map must contain at least one dissipative repeller. These attractors and repellers are
different from $C$ (because they are dissipative) and stay at a finite distance from $C$ (because they are
compact sets). Thus, we may take smaller absorbing domains around $C$ and obtain one more dissipative attractor even closer to $C$, and a dissipative repeller as well. This procedure can be repeated infinitely many times, which gives us the following result.\\

\noindent{\bf Theorem 1.} {\em Every reversible core contains a limit of an infinite sequence of dissipative attractors and an infinite sequence of dissipative repellers.}\\

As we have attractors, repellers, and a kind of conservative object (the reversible core) unseparated from each other, we speak about mixed dynamics near the core. The numerical detection of a reversible core can be based on the following\\

\noindent{\bf Theorem 2.} {\em If a CRH-attractor has a non-empty intersection with a CRH-repeller, they must coincide and form a reversible core.}\\

{\em Proof.} Since the attractor $A$ and the repeller $R$ are both chain-transitive and have a non-empty intersection, their union is also chain-transitive. Therefore, every point of $A\cup R$
is attainable by $\varepsilon$-orbits that start in $A$ for every $\varepsilon>0$, which means $A=A\cup R$ (because $A$, by
definition, is a stable set, meaning that $\varepsilon$-orbits that
start in $A$ cannot get far from $A$). Similarly, $R=A\cup R$, so $A$ and $R$ coincide. Now, since $A$ is a CRH-repeller, no $\varepsilon$-orbit that
starts at a bounded away from zero distance from $A$ can get close to $A$ if $\varepsilon>0$ is sufficiently small. Thus, $A=R$ is a reversible core, and
both forward and backward $\varepsilon$-orbits never leave its small neighbourhood. $\square$\\

By Theorem 2, if we have a reversible core, then in numerical simulations we would see that the attractor and repeller occupy roughly the same region in the phase space.
However, we will not see that the attractor coincides exactly with a repeller:
as numerics add some small noise, it makes sense to expect that a numerical forward orbit will shadow some absorbing domain around the core,
while the backward numerical orbit will shadow some absorbing domain for the inverse map; such domains can not completely coincide.
Thus, a numerical indication of the mixed dynamics would be a numerically obtained attractor which would not coincide
with a numerically obtained repeller but would have a sizable intersection with it and the difference between these two sets
would appear small.
\begin{figure}[h!]
\centering
\includegraphics[width=12cm]{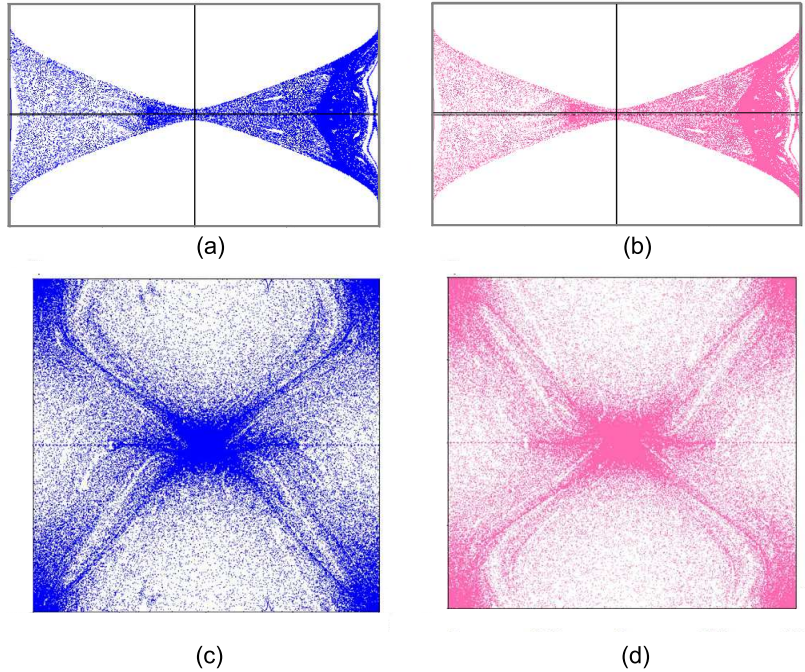}
\caption{{\footnotesize Attractor (a) and repeller (b) for the the Suslov model (the pictures are taken from \cite{K16}). Attractor (c) and repeller (d) for the the Pikovsky-Topaj model
(the pictures are taken from \cite{GGKT17}). }}\label{suslov}
\end{figure}

\textbf{Remark.}
The question therefore arises:
do Figs.~1b and 1d indicate the presence of a large reversible core? Note that
the systems whose phase portraits are depicted in Fig.1 are not covered by the theory developed in the next Sections, as  they cannot have generic elliptic points (the dimension of $\mbox{Fix}\;(g)$ is less than half
dimension of the phase space -- for example, in the model of ``rubber ellipsoid'' (Fig.~1c,d), the involution $g$ that takes attractor to repeller is
the central symmetry, so it is orientation-preserving and $\mbox{Fix}\;(g)$ is just one point). On the other hand, there are various examples of reversible systems having symmetric generic elliptic points, see e.g. \cite{PikTop02,GGKT17,K13,K16}. In particular, in the Suslov model \cite{K16} and in the Pikovsky-Topaj model \cite{PikTop02}, the involution $g$ of the corresponding two-dimensional Poincar\'e map   is $x\to x,\; y\to -y$ and, thus, the attractor and repeller are always symmetric with respect to the axes $y=0$, see Fig.~\ref{suslov}. The question of whether the ``large'' intersection of the attractor and repeller observed in Figs. 1 and 5 is a reversible core, is it related to elliptic orbits in the cases of Fig. 5, and with which it is related in the cases of Fig.1 remains open.

\section{Generic reversible cores in two-dimensional reversible maps}
\label{el2}

We now show that reversible cores exist for a large class of dynamical systems. We restrict ourselves to reversible
maps, as they are known to provide examples of mixed dynamics in abundance, as we mentioned in the Introduction.
Moreover, in this Section we consider only two-dimensional maps. Thus, let $f$ be a $C^r$-diffeomorphism
($r=1,\dots,\infty$) of a two-dimensional orientable manifold and assume that $f$ is reversible, i.e.,
\begin{equation}\label{rdf2}
f^{-1}=g\circ f\circ g,
\end{equation}
where $g$ is a $C^r$-smooth involution (a map such that $g\circ g=id$).

A periodic orbit $\{x_0,\dots,x_m\}$ of $f$ is called symmetric if it is invariant with respect to $g$; namely, $gx_0=f^jx_0$ for some $j\leq m$ (then, by (\ref{rdf2}), $gfx_0=f^{j-1}x_0$, and so on). It is easy to see that for a symmetric
periodic orbit at least one of its points is either a fixed point of $h=g$, or a fixed point of $h=f\circ g$. Let $x_0$ be such point -- we call it a symmetric periodic point.
By (\ref{rdf2}), $h$ is an involution and
\begin{equation}\label{rdfh2}
T^{-1}=h\circ T\circ h,
\end{equation}
where $T=f^m$ is the first-return map near $x_0$ (so $Tx_0=x_0$).

For a symmetric periodic orbit,
\begin{itemize}
\item
if $\lambda$ is a multiplier, then $\lambda^{-1}$ is also a multiplier.
\end{itemize}

Indeed, let $x_0$ be a symmetric periodic point, i.e. $T x_0=x_0$ and $h x_0= x_0$.
Denote as $A=T'$ the derivative of $T$ at $x_0$. By Bochner theorem \cite{Bo}, we may always choose coordinates near $x_0$ such that the involution  $h$ is linear. By (\ref{rdfh2}), we have
\begin{equation}\label{adfh2}
A^{-1}=h\circ A\circ h.
\end{equation}
If $A{\bf e}=\lambda {\bf e}$, i.e. $\bf e$ is an eigenvector of $A$ with the eigenvalue (the multiplier) $\lambda$,
then it follows from (\ref{adfh2}) that $A^{-1} h{\bf e}= \lambda h{\bf e}$, i.e. $h{\bf e}$ is also an eigenvector of $A$ with the eigenvalue $\lambda^{-1}$.

Note that if $\lambda^2\neq 1$, i.e., the multipliers $\lambda$ and $\lambda^{-1}$ are different, then the eigenvectors
$\bf e$ and $h{\bf e}$ are not collinear. Thus, the involution $h$ interchanges a pair of non-collinear vectors, which means
that $h$ reverses the orientation. We will further assume that the map $f$ is orientation-preserving, so
the original involution $g$ must be orientation-reversing in this case (recall that $g=h$ or $g=f^{-1}\circ R$). This
will be our standing assumption from now on.

Since the linearization matrix $A$ is real, if $\lambda$ is its eigenvalue, than the complex conjugate $\lambda^*$ must also be an eigenvalue. So, if the symmetric periodic orbit
has a complex (not real) multiplier $\lambda$, then $\lambda^*=\lambda^{-1}$, i.e., both multipliers must lie on
the unit circle. This means that there exists $\omega\in(0,\pi)$ such that the multipliers of the periodic orbit are
$e^{\pm i\omega}$. We call the symmetric periodic orbit {\em elliptic} in this case.

It is well-known \cite{Sevr} that a symmetric elliptic periodic orbit of a two-dimensional reversible $C^r$-diffeomorphism remains elliptic
at $C^r$-small perturbations which keep the map reversible. So, systems with elliptic orbits form an open subset in the space of
$g$-reversible $C^r$-diffeomorphisms. Empirical evidence suggests that this open set should be quite large, cf. \cite{GGKT17}. In particular, elliptic orbits are born \cite{LS04} at bifurcations of reversible maps with heteroclinic cycles described in the Introduction and Section \ref{absnew}. Therefore, systems with elliptic periodic orbits form a dense (and open) subset of the absolute Newhouse domain in the space ${\cal R}^r_g$ of two-dimensional $g$-reversible $C^r$-diffeomorphisms. On the other hand, bifurcations near elliptic orbits lead to creation of ultimately-wild hyperbolic sets \cite{GLRT14}, which means that the $C^r$-closure of the set of systems from ${\cal R}^r_g$ with elliptic orbits coincides with the $C^r$-closure of the absolute Newhouse domain in ${\cal R}^r_g$.

A natural conjecture (a similar conjecture
for area-preserving maps can be found e.g. in \cite{N77}) would be that if the map has a chaotic attractor but
is not uniformly-hyperbolic, then homoclinic tangencies or cycles with heteroclinic tangencies can be created by
$C^r$-small perturbations. No mathematical technique is currently available for a proof of such statement. Still it is reasonable to conjecture that if a chaotic attractor is observed in a given $g$-reversible map,
then most probably the attractor is not uniformly-hyperbolic and contains a wild hyperbolic set (unless we deal with an Anosov map on a torus), and if the attractor intersects its own image by the involution $g$, then the wild-hyperbolic set
is ultimately-wild (i.e., it contains periodic orbits with Jacobians both larger and smaller than $1$). Thus, such map can be suspected to be in the absolute Newhouse domain and, in the case of orientation-reversing involution $g$, this means
that elliptic periodic orbits should be expected.

By \cite{Sevr} most of the neighborhood of a generic\footnote{The genericity conditions in this case include the absence of strong
resonances, i.e. $\omega\neq \pi/2, 2\pi/3$, and the non-vanishing of one of the Birkhoff coefficients in the case of irrational
$\omega/\pi$ or, in the case of a rational $\omega/\pi$, one of the first Birkhoff coefficients. In the conservative
case the generic elliptic point is KAM-stable \cite{Sevr}. In the reversible case, the KAM-stability, of course, holds, but,
as Theorem 3 shows, elliptic points of a generic non-conservative reversible map are stable in a stronger sense.} elliptic point of a reversible map is occupied by invariant KAM-curves, like in the
conservative case. However, as it was shown in \cite{GLRT14}, a generic elliptic point of a two-dimensional non-conservative reversible map is also a limit of infinite sequences of periodic attractors and periodic repellers (that are born in the resonant zones). Here, we strengthen this result and prove the following\\

{\bf Theorem 3.} {\em All symmetric elliptic periodic orbits of a $C^r$-generic two-dimensional $g$-reversible map are reversible cores.}\\

{\em Proof.} Since a periodic orbit is a chain-transitive set, it suffices to prove that each elliptic orbit is (generically) surrounded by a sequence of nested absorbing domains both for the map $f$ itself and its inverse. We start with
proving that given an elliptic orbit $P$ and any its open neighborhood $U$ one can make an arbitrarily small $C^r$-perturbation of the map (within the class ${\cal R}^r_g$ of reversible systems) such that the perturbed map $f$
and its inverse $f^{-1}$ would have, each, an absorbing domain lying inside $U$ and containing $P$.

Let $x_0$ be a symmetric point on $P$ and $T$ be the first-return map near $x_0$. We can introduce a complex
coordinate $z$ near $x_0$ such that the map $T$ will be the period map of the time-reversible non-autonomous flow \cite{Lamb} defined by the differential equation
\begin{equation}\label{eq1}
\dot z = i\omega z + F(z,z^*,t)
\end{equation}
where the $C^r$-function $F$ is 1-periodic , i.e., $F(z,z^*,t)=F(z,z^*,t+1)$, and it vanishes at $z=0$ along with its derivative with respect to $z$ and $z^*$, so $z=0$ is the elliptic fixed point of $T$. The reversibility means here that
$F(z,z^*,t)=F^*(z^*,z,-t)$, i.e., this system is invariant with respect to the transformation
$\{t\to -t, \; z \leftrightarrow z^*\}$, so the involution $h$ from (\ref{rdfh2}) is here the complex conjugation
operation $\{z \to z^*\}$.

One can add an arbitrarily small perturbation to (\ref{eq1}) such that $\omega/\pi$ would become irrational. Then
the normal form theory for reversible maps \cite{Lamb} will ensure the existence of a $C^r$ coordinate transformation
which brings the map $T$ to the form
\begin{equation}\label{prd}
T=R_\omega \circ {\cal T}
\end{equation}
where $R_\omega: (z,z^*) \mapsto (e^{i\omega}z, e^{-i\omega}z^*)$ is the rotation to the angle $\omega$, and
$\cal T$ is the time-1 map by the flow of
\begin{equation}\label{eq2}
\dot z = i \Omega(|z|^2)z + o(|z|^r),
\end{equation}
where $\Omega$ is a real polynomial such that $\Omega(0)=0$; all time-dependent terms are now in the $o(|z|^r)$-term
(which vanishes at $z=0$ along with all derivatives with respect to $z$ and $z^*$ up to the order $r$). The next step is to
remove the $o(|z|^r)$-term in (\ref{eq2}), i.e., add a $C^r$-small perturbation after which the map $\cal T$
will, in a sufficiently small neighborhood of zero, coincide with the time-1 map of the reversible autonomous flow given by
\begin{equation}\label{eq3}
\dot z = i \Omega(|z|^2)z .
\end{equation}
By an additional small perturbation we can always achieve that $\Omega_1=\Omega'(0)\neq 0$.

By a $C^r$-small perturbation of the map $f$ within the class of reversible maps we can, while keeping the map $T$
in the form (\ref{prd}), make $\omega=2\pi \frac{p}{q}$ where $p$ and $q$ are co-prime integers
(we will also assume that $q\geq r$). We may also change equation (\ref{eq3}) such that in a small neighborhood of $z=0$
\begin{equation}\label{eq4}
\dot z = -i\mu z + i \Omega(|z|^2)z + i \delta (z^*)^{q-1} + iBz^{q+1}+iCz(z^*)^q,
\end{equation}
where $\delta$ and $\mu$ are small real parameters and $B$ and $C$ are some real constants. Note that all
coefficients in the right-hand side of (\ref{eq4}) are pure imaginary, so the equation is time-reversible, as it should.

Note that equation (\ref{eq4}) is the truncated normal form for the bifurcations of a symmetric elliptic point with
a degenerate resonance, corresponding to the multiplier $\lambda=e^{2\pi i \frac{p}{q}}$ and the coefficient of the
first non-trivial resonant term $(z^*)^{q-1}$ vanishing at the moment of bifurcation. Without the $O(|z|^{q+1})$-terms,
equation (\ref{eq4}) is Hamiltonian. However, at non-zero $B$ and $C$, the conservativity is broken -- the equation may
have e.g. asymptotically stable and unstable equilibria, as we see below.

The equation is symmetric with
respect to rotation to the angle $2\pi \frac{p}{q}$. This means that the maps $R_{2\pi p/q}$ and $\cal T$ in (\ref{prd})
commute, so the map $T^q$ is the time-~$q$ map of (\ref{eq4}). In particular, every equilibrium state of (\ref{eq4}) is
a period-$q$ point of $T$; asymptotically stable equilibria correspond to asymptotically stable periodic points. In this way,
based on the analysis of stability of equilibrium states of (\ref{eq4}), it was shown in \cite{GLRT14} that stable periodic orbits can be born at bifurcations of elliptic periodic points in reversible maps. Here, we need a more detailed investigation of
the dynamics of equation (\ref{eq4}).

To this aim, introduce polar coordinates: $z=\sqrt{\rho} e^{i\phi/q}$. The equation takes the form
\begin{equation}\label{rph} \begin{array}{l}\displaystyle
\dot \rho = 2\rho^{\frac{q}{2}} (\delta - (B-C)\rho) \sin\phi, \\ \displaystyle
\dot\phi = q(\Omega(\rho) - \mu) + q\rho^{\frac{q-2}{2}} (\delta + (B+C)\rho) \cos\phi.\end{array}
\end{equation}
Note that this system is invariant with respect to $\{t\to -t, \;\phi\to-\phi\}$.

\begin{figure}[h!]
\centering
\includegraphics[width=14cm]{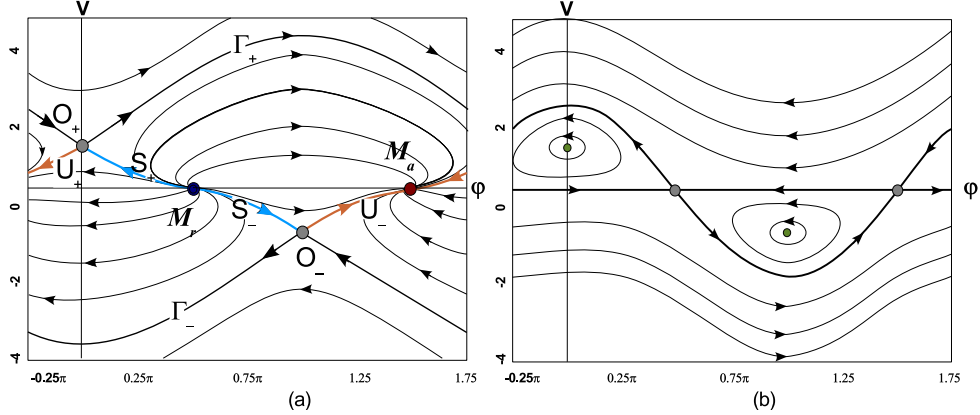}
\caption{{\footnotesize The phase portrait of system (\ref{r0ph0}) on the cylinder for (a) $D=0$, $\beta=1$ and (b) $D=0$, $\beta=-1$.}}
\label{figD0PP}
\end{figure}

Assume $B\neq 0$ and $B\neq C$. Choose sufficiently small
$\rho_0>0$, put
\begin{equation}\label{mrho}
\mu=\Omega(\rho_0),
\end{equation}
and consider the behavior of the system for $\rho$ close to $\rho_0$. We do this by scaling
\begin{equation}\label{rhrho}
\rho=\rho_0 - \rho_0^{q/2} \frac{2B}{\Omega_1}V,
\end{equation}
where the range of values of $V$ can be as large as we want if $\rho_0$ is small enough (the coefficient $\Omega_1\neq 0$ equals to $\Omega'(0)$). We also scale the small parameter
\begin{equation}\label{deltarh}
\delta=(B-C) \rho_0 +\frac{2B}{\Omega_1} (B-C) \rho_0^{q/2} D,
\end{equation}
where the rescaled parameter $D$ is no longer small
and can take arbitrary finite values, and introduce the
new time $s=2(C-B) \rho_0^{q/2}\;t$. The system will take the form
\begin{equation}\label{rph0} \begin{array}{l}\displaystyle
\dot V =  (D + V) \sin\phi + O(\rho_0^{\frac{q-2}{2}}), \\ \displaystyle
\dot\phi =\frac{B q}{B-C}(V  - \cos\phi)+O(\rho_0^{\frac{q-2}{2}}).\end{array}
\end{equation}
The limit of this rescaled system as $\rho_0\to 0$ is
\begin{equation}\label{r0ph0}
\dot V =  (D + V) \sin\phi , \qquad \dot\phi =\beta (V  - \cos\phi),
\end{equation}
where $\beta=q B/(B-C)$. This is a time-reversible system on the cylinder parameterized by $(V,\phi)$.
Importantly, this system can be solved. Thus, it is easy to see
that the phase curves of this system satisfy
$$cos\phi + D = (D+V) \frac{\beta}{\beta-1}+ K |D+V|^{\beta}$$
with indefinite constants $K$. At $\beta=1$ this formula should be replaced by
$$cos\phi + D = K (D+V) - (D+V) \ln |D+V|.$$
With these formulas, one can construct the phase portrait of system (\ref{r0ph0}) on the cylinder for different values
of $D$ and $\beta$, as shown in Fig.~\ref{figD0PP}.

In particular, for $\beta>0$ and $|D|<1$, see Fig.~\ref{figD0PP}a, this system has two symmetric (with respect to the involution $\phi\to -\phi$)
saddle equilibria $O_+(1,0)$ and $O_-(-1,\pi)$ and two asymmetric equilibria $M_{a,r}=(-D,\phi_{a,r})$ where
$\cos\phi_{a,r}=-D$, the equilibrium $M_a$ (with $\sin\phi_a<0$) is asymptotically stable and $M_r$ (with $\sin\phi_r>0$) is asymptotically unstable. Two of the separatrices of $O_+$ coincide and form a homoclinic loop $\Gamma_+$, while of the two other separatrices the unstable separatrix $U_+$ tends to $M_a$ as $t\to+\infty$ and the stable separatrix $S_+$ tends to $M_r$ as $t\to-\infty$. The same is true for $O_-$: two of its separatrices form a homoclinic loop $\Gamma_-$ and the two other separatrices tend one ($U_-$) to $M_a$ as $t\to+\infty$
and the other ($S_-$) to $M_r$ as $t\to-\infty$. In the invariant annulus bounded by $\Gamma_+$ and $\Gamma_-$ all the
orbits in its interior, except for the repeller $M_r$ and the separatrices $S_\pm$, tend to $M_a$ as $t\to+\infty$, while
all the orbits except for the attractor $M_a$ and the separatrices $U_\pm$ tend to $M_r$ as $t\to-\infty$.

\begin{figure}[h!]
\centering
\includegraphics[width=14cm]{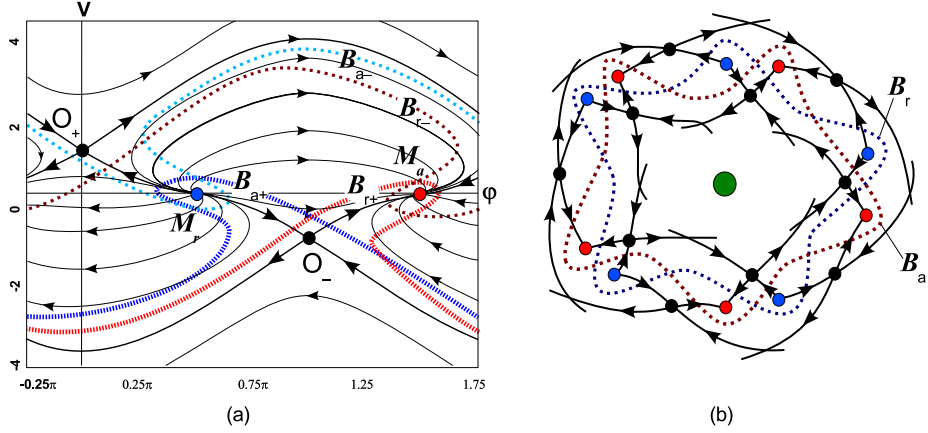}
\caption{{\footnotesize (a) Boundaries of absorbing domains for system (\ref{r0ph0}) and its inverse (for $t\to -t$). The domains with boundaries $B_{a+}$ and $B_{r+}$ contain, respectively, the attractor $M_a$ and the repeller $M_r$ and the upper part of the cylinder.
The domains with boundaries $B_{a-}$ and $B_{r-}$ contain, respectively, $M_a$ and $M_r$ and the lower part of the cylinder.
(b) A pair of absorbing domains with boundaries $B_a$ and $B_r$ around the point $z=0$ (for system (\ref{eq4}) and for its time reversal).}}
\label{figD0At3}
\end{figure}

It follows
that if we remove from the phase cylinder a small neighborhood of $\Gamma_+\cup S_+\cup M_r$, then the connected
component that contains $M_a$ is an absorbing domain. Note that it contains the part of the cylinder corresponding to $V\to-\infty$. Similarly, one obtains an absorbing domain corresponding to $V\to +\infty$ -- this is the connected component of the cylinder minus a small neighborhood of $\Gamma_-\cup S_-\cup M_r$, which contains $M_a$ (see Fig.~\ref{figD0At3}a). In the same way, by removing a small neighborhood of $\Gamma_+\cup U_+\cup M_a$ or a small neighborhood of $\Gamma_-\cup U_-\cup M_a$ from the cylinder and taking the connected component that contains the repeller $M_r$, we obtain a pair of absorbing domains (with the boundaries $B_{r+}$ and $B_{r-}$, see Fig.~\ref{figD0At3}a)for the system obtained from (\ref{r0ph0}) by the time reversal.

The absorbing domains do not disappear at small perturbations of the system, so they persist
for system (\ref{rph0}) for all small $\rho_0$. System (\ref{rph0}) is obtained from (\ref{eq4}) by a rescaling of coordinates. In the non-rescaled coordinates, $z=0$ corresponds either to very large positive $V$ or to very large negative $V$, so in any case we obtain that, for appropriately chosen
values of $\mu$, $\delta$, $B$, and $C$, both
system (\ref{eq4}) and the system obtained from it by the time reversal have, each, an absorbing domain contains the equilibrium at $z=0$ . Note that $\mu$ and $\delta$ can be made as small as we want by taking $\rho_0$ small (see (\ref{mrho}),(\ref{deltarh})), as required. Note also that the absorbing domains
containing points with large negative values of $V$ do not extend beyond $\Gamma_+$ and the absorbing domains
containing points with large positive values of $V$ do not extend beyond $\Gamma_-$. As these boundaries
correspond to bounded values of $V$, it follows that in the non-rescaled coordinates $z$ these boundaries correspond
to $|z|$ close to $\rho_0$, i.e., the pair of absorbing domains around the point $z=0$ (for system (\ref{eq4}) and for its time reversal) lie entirely in an $O(\rho_0)$-neighborhood of this point, see Fig.~\ref{figD0PP}b.

The time-$q$ map of system (\ref{eq4}) is the map $T^q$. Thus, we have constructed a $C^r$-small perturbation of the map $f$ such that the map $T^q=f^{qm}$ (where $m$ is the period of the periodic orbit $P$ under consideration) and the map $T^{-q}$ have, each,
an absorbing domain ${\cal D}_a$ and, respectively, ${\cal D}_r$, in a small (as small as we want) neighborhood of some point on $P$. Obviously, a small
open neighborhood of the closure of $\bigcup_{i=0}^{qm-1} f^i  {\cal D}_a$ is an absorbing domain for the map $f$, which
contains the whole orbit $P$ and is contained in a small neighborhood of $P$. A small
open neighborhood of the closure of $\bigcup_{i=0}^{qm-1} f^{-i}  {\cal D}_r$ is an absorbing domain for the map $f^{-1}$, it contains the whole orbit $P$ and is contained in a small neighborhood of $P$.

We can now finish the proof of the theorem. First, we recall the well-known fact that if a two-dimensional $g$-reversible map with an orientation-reversing involution $g$ has a symmetric orbit $L$ of period $m$ with the multipliers $\lambda_1=\lambda_2=1$ or $\lambda_1=\lambda_2=-1$, then by an arbitrarily small $C^r$-perturbation of the map
one can achieve that in a small neighborhood of $L$ all orbits of the same period $m$ will be either elliptic or hyperbolic
(i.e., $|\lambda_{1,2}|\neq 1$). Now, take a countable base of balls $U_s$, $s=1,\dots,\infty$ (so every open set is a union of some sequence of the balls $U_s$). Chose one of these balls, and take an integer $m\geq 1$. As we just mentioned, by
an arbitrarily small perturbation of the map $f$ within the class ${\cal R}^r_g$ one can achieve that all symmetric orbits of period $m$ or less that intersect the chosen ball $U_s$ are either elliptic or hyperbolic. There are no other orbits of the same or smaller period in a neighborhood of an elliptic or hyperbolic periodic orbit of period $m$, nor
such orbits can be born at a $C^r$-small perturbation, so
the number of points of symmetric elliptic orbits of period $\leq m$ in $U_s$ is finite and this property holds for an
open and dense set for maps from ${\cal R}^r_g$. As we proved above, by an additional $C^r$-small perturbation
we can create a pair of absorbing domains, one for the map $f$ and the other for the map $f^{-1}$, around each of the elliptic orbits of period $\leq m$ that intersect $U_s$, and these absorbing domains lie inside
$\bigcup_{i=-mr}^{mr} f^i(U_s)$. Absorbing domains persist under small perturbations of the map, so the set ${\cal E}_{s,m}$
with this property is open and dense in ${\cal R}^r_g$. The intersection ${\cal E}^*$ of the sets ${\cal E}_{s,m}$ over all $U_s$ and all integer $m$ is a countable intersection of open and dense sets, so it is a residual set and every map from this intersection is, by definition, $C^r$-generic. By construction, for every map from ${\cal E}^*$, for every elliptic orbit of it,
in any neighborhood of this orbit there exists a pair of absorbing domains, one for the map itself, one for the inverse map, such that both domains contain this orbit. Thus every elliptic orbit of every map from ${\cal E}^*$ is, simultaneously
a CRH-attractor and a CRH-repeller, i.e., it is a reversible core. $\square$.

\section{Universal dynamics near elliptic orbits in reversible systems}
\label{eln}

The richness of dynamics near a reversible core does not need to be exhausted by simple periodic attractors and repellers.
In particular, as we show in this Section, near elliptic orbits of reversible maps, the dynamics can be as complicated and diverse as it only possible for the given dimension of the phase space. We will not restrict our consideration to two-dimensional maps here, so we start with reminding the classification of periodic orbits of $n$-dimensional
reversible maps.

Let ${\cal M}$ be an $n$-dimensional manifold. Let $g:{\cal M}\rightarrow {\cal M}$ be an involution, so $g\circ g=id$.
A diffeomorphism $f:{\cal M}\rightarrow {\cal M}$ is called a reversible map if it is conjugate by $g$ to its own inverse.
A periodic orbit of $f$ is called symmetric if it is invariant with respect to $g$; at least one point $x_0$
of a symmetric periodic orbit is a fixed point of an involution $h$ where $h=g$ or $h=f\circ g$.
Let $T$ be the first-return map near $x_0$, so $T x_0=x_0$ and $h x_0= x_0$. By reversibility
\begin{equation}\label{rdfh}
T^{-1}=h\circ T\circ h.
\end{equation}
Denote as $A$ the derivative of $T$ at $x_0$. We may always choose $C^r$-coordinates near $x_0$ such that the involution $h$ is linear \cite{Bo}. By (\ref{rdfh}), we have
\begin{equation}\label{adfh}
A^{-1}=h\circ A\circ h.
\end{equation}
If $A{\bf e}=\lambda {\bf e}$, i.e. $\bf e$ is an eigenvector of $A$ with the eigenvalue (the multiplier) $\lambda$,
then it follows from (\ref{adfh}) that $A^{-1} h{\bf e}= \lambda h{\bf e}$, i.e. $h{\bf e}$ is also an eigenvector of $A$
with the eigenvalue $\lambda^{-1}$.

Note that if an eigenvector of $A$ corresponds to a multiplier which is not $\pm 1$ (i.e. $\lambda\neq \lambda^{-1}$), then this eigenvector
is not an eigenvector of the involution $h$. On the other side, every eigenvector of $A$ which corresponds to $\lambda=\pm1$ is a linear combination of
eigenvectors which are, at the same time, eigenvectors of $h$.\footnote{If $A{\bf e}=\lambda {\bf e}$ with $\lambda^2=1$ and $h{\bf e}\neq \pm {\bf e}$, then ${\bf e}+h{\bf e}$ and
${\bf e}-h{\bf e}$ are linearly independent eigenvectors of $A$. Indeed, since $h{\bf e}$ is an eigenvector of $A$ with the eigenvalue
$\lambda^{-1}=\lambda$, we have $A({\bf e} T- h{\bf e}) = A{\bf e} T- Ah{\bf e} = \lambda ({\bf e} T- h{\bf e})$.}
One can also show that if $A$ has an eigenvector ${\bf e}$ with the eigenvalue $\pm 1$ and this eigenvector is not an eigenvector of $h$,
then one can, without destroying the reversibility, add an arbitrarily small perturbation to the map $T$ such that, for the perturbed $T$, we would have
$Ay=\tilde\lambda {\bf e}$, $A \cdot (h{\bf e})=\tilde\lambda^{-1}\; h{\bf e}$ for some $\tilde\lambda$ different from $\pm 1$. This means that for a generic symmetric periodic orbit,
exactly those multipliers that correspond to the eigenvectors of $A$ which are not eigenvectors of $h$ are different from $\pm 1$.

One can also show that the linearization matrix $A$ for a generic symmetric periodic orbit has no Jordan blocks. Thus, the invariant subspace $I_+$
of the matrix $A$ which corresponds to the eigenvalue $+1$ is spanned by eigenvectors of $A$, and either all of them satisfy $h{\bf e}={\bf e}$ or all of them satisfy $h{\bf e}=-{\bf e}$.
The same is true for the invariant subspace $I_-$ which corresponds to the multiplier $-1$. However, we can always, if necessary, replace the involution $h$
by the involution $T\circ h$ (identity (\ref{rdfh}) would not change) and achieve that the involution is identity on $I_-$.
After this choice is made, it can be shown (see remark after formula (\ref{pnfff00}) in the proof of Theorem 5 below) that if $h=-id$ on $I_+$, then the symmetric
periodic point can be made to disappear by an arbitrarily small perturbation of the map. Therefore, for the generic
symmetric periodic orbit we have
\begin{equation}\label{pnhid}
h|_{I_+\oplus I_-}=id.
\end{equation}

It can be shown that condition (\ref{pnhid}) is necessary and sufficient for the symmetric periodic orbit to persist at
small smooth perturbations which preserve the reversibility. Below, such periodic orbits will be called {\em regular}.

If the matrix $A$ has some hyperbolic eigenvalues, i.e. the multipliers not on the unit circle, then the map $T$ near the periodic point $x_0$ has an invariant center manifold $W^c$, which is an intersection of the center-stable manifold $W^{cs}$ and the center-unstable
manifold $W^{cu}$. The invariant manifold $W^{cs}$ is tangent to the invariant space of $A$ which corresponds to all the multipliers smaller or equal to $1$ in the absolute value (i.e., the multipliers on the unit circle or inside it), while the invariant manifold $W^{cu}$ is tangent to the invariant space of $A$ which corresponds to
all the multipliers whose absolute value is greater or equal to $1$ (the multipliers on the unit circle and outside of it).
We can choose $W^{cu}=h(W^{cs})$, then the invariant center manifold $W^c=W^{cs}\cap h(W^{cs})$
will be $h$-invariant. The center manifold $W^c$ is smooth (of class $C^r$ for any finite $r$) and persists at $C^r$-small perturbations. For every map close to $f$, all the orbits that never leave a small neighborhood of the periodic orbit under consideration must belong to $W^c$. Dynamics transverse to $W^c$ is hyperbolic (trivial) -- the orbits not from $W^c$ either leave a small neighborhood of the periodic orbit both at forward and backward iterations of $f$, or they lie in
$W^{cs}$ and exponentially approach $W^c$ at forward iterations of $f$ (so they leave a small neighborhood of the periodic orbit at backward iterations), or they lie in $W^{cu}$ and exponentially approach $W^c$ at backward iterations of $f$ (and leave a small neighborhood of the periodic orbit at forward iterations). Dynamics in $W^c$ is, generically,
very non-trivial. In particular, we have the following result.\\

{\bf Theorem 4.} {\em Consider a regular symmetric periodic orbit of a $g$-reversible $C^r$-smooth map $f$.
Assume the orbit has at least one
pair of complex multipliers on the unit circle, i.e., $\lambda=e^{\pm i\omega}$ where $\omega\in(0,\pi)$.
Then for a $C^r$-generic $g$-reversible map sufficiently close to $f$ the periodic orbit is a limit of an infinite sequence of uniformly hyperbolic attractors and uniformly hyperbolic repellers of all topological types
possible for a smooth map of a $d$-dimensional disc, where $d$ is the dimension of the center manifold for this orbit.}\\

This theorem is an immediate consequence of a more general statement (Theorem 5) that employs the notion of a {\em universal} map from \cite{T03,T15}. Given a $C^r$-smooth diffeomorphism $F$ of an $n$-dimensional manifold $\cal M$,
we consider the set of the so-called renormalized iterations of $F$, defined as follows. Take a unit ball ${\cal B}_n\subset R^n$. Let $\psi$ be a $C^r$-map $R^n\to {\cal M}$, which is a diffeomorphism between $R^n$ and its image
$\psi(R^n)$. Take the ball $\psi({\cal B}_n)$ and suppose that its image $F^k\circ \psi({\cal B}_n)$ lies inside
$\psi(R^n)$ for some positive integer $k$. Then the $C^r$-diffeomorphism $F_{k,\psi}: {\cal B}_n\to R^n$, defined
as $F_{k,\psi}=\psi^{-1}\circ F^k\circ \psi|_{{\cal B}_n}$ is a renormalized iteration of $F$.
\begin{itemize}
\item
The map $F$ is called
$d$-{\em universal} if the set of its renormalized iterations is $C^r$-dense in the set of all orientation preserving $C^r$-diffeomorphisms from ${\cal B}_n$ to $R^n$.
\end{itemize}

By the definition, iterations of any $d$-universal map approximate arbitrarily well all dynamics possible in a $d$-dimensional ball. Therefore, every $C^r$-robust phenomenon occurring in any $d$-dimensional diffeomorphism is also present in
each $d$-universal map. In particular, every $d$-universal map has, simultaneously, uniformly-hyperbolic attractors and
repellers of all topological types possible in the $d$-dimensional ball. Thus, Theorem 4 is a direct consequence
of the following result, which we obtain by using the theory from \cite{T15} about smooth perturbations of identity.\\

{\bf Theorem 5.} {\em Consider a regular symmetric periodic orbit of a $g$-reversible $C^r$-smooth map $f$.
Assume the orbit has at least one
pair of complex multipliers on the unit circle. Then any $C^r$-generic $g$-reversible map sufficiently close to $f$, restricted to the local center manifold in an arbitrarily small neighborhood of the periodic orbit under consideration, is $d$-universal, where $d$ is the dimension of the center manifold.}\\

{\em Proof.} Take a regular symmetric periodic point $x_0$  and consider a small piece of the local center manifold $W^c$ that contains $x_0$. We will show below that by an arbitrary small perturbation of the map within the class of $C^r$-smooth $g$-reversible maps one can create a {\em periodic spot} in $W^c$, i.e., an open $d$-dimensional ball in $W^c$
for which all points are periodic {\em with the same period}. After that is done, we will just need an easy adaptation of
the result of \cite{T15} to the reversible case in order to make the map universal.

Thus, we have a periodic orbit with the multipliers $e^{\pm i\omega_1}, \dots, e^{\pm i\omega_m}$ ($0<\omega_j<\pi$, $m\geq 1$), $s$ multipliers equal to $+1$, $k$ multipliers equal to $-1$, so $m+k+s=d$ multipliers lie on the unit circle, and the rest $(n-d)$ of the multipliers not on the unit circle. As it was explained before, we may generically assume
that there are no Jordan blocks in the linearization matrix $A$ of the period map $T$ (for multipliers $\pm 1$ this is due to (\ref{pnhid});
for the other multipliers we assume that they are all simple). Also, by a small perturbation, we can make the values of $\omega_j$, $j=1,\dots,m$, rationally independent with each other and $\pi$.

Let $\sigma$ be the linear map that changes sign of the coordinates in the invariant subspace of the matrix $A$ that
corresponds to all real negative multipliers. Then, the normal form for the map $T$ near $x_0$ equals \cite{Lamb} to $\sigma$ times the time-$1$ map by the flow of a system of differential equations of the form
\begin{equation}\label{pnfff}
\begin{array}{l}\dot z_j=i \Omega_j(Z,u,v)z_j + O(y)+o(\|x-x_0\|^r) \qquad (j=1,\dots, m)\\
\dot u= F(Z,u,v) + O(y)+o(\|x-x_0\|^r), \\ 
\dot v =G(Z,u,v) +O(y)+o(\|x-x_0\|^r), \\
\dot y=H(Z,u,v,y)y +o(\|x-x_0\|^r),\end{array}
\end{equation}
where $x=x_0+(z,u,v,y)$, the variables $z_j$ are complex (projections to the eigenspaces corresponding to the multipliers
$e^{\pm i\omega_j}$), $Z=(|z_1|^2,\dots, |z_s|^2)$, $u\in R^s$ (these are projections to the eigenspace corresponding to the multipliers $+1$), $v\in R^k$ (these are projections to the eigenspace corresponding to the multipliers $-1$), and
$y\in R^{n-d}$ (projections to the eigenspace corresponding to the multipliers not on the unit circle);
the functions $F$ and $G$ have zero linear parts, $\Omega_j(0,0,0)=\omega_j$, and
$H(0,0,0,0)$ is a matrix with eigenvalues outside of the imaginary axis. The $o(\|x-x_0\|^r)$ terms can be time-dependent
(anti-periodic in the sense that they do not change after multiplication to $\sigma$ and adding $1$ to time); they vanish at zero along will the derivatives up to the order $r$. The $O(y)$ terms stand for the functions vanishing at $y=0$.
The normal form system is also $\sigma$-equivariant, e.g.
$\Omega_j$ and $F$ are even functions of $v$, and $G$ is an odd function of $v$.
The same normalizing transformation brings the involution $h$ from (\ref{rdfh}),(\ref{adfh}) to the form $h(z,u,v,y)=(z^*,u,v, \hat h y)$ or $h(z,u,v,y)=(z^*,-u,v, \hat h y)$, where $z^*$ is complex conjugate to $z$ and $\hat h$ is some linear involution in the $y$-space.

By adding a $C^r$-small perturbation to $f$ one can make the $o(\|x-x_0\|^r)$-terms in (\ref{pnfff}) vanish
in a sufficiently small neighborhood of $x_0$. Then, the map $T$ near $x_0$ will be equal to $\sigma$ times the time-$1$ map by the flow of the autonomous system
\begin{equation}\label{pnfff0}
\begin{array}{l}\dot z_j=i \Omega_j(Z,u,v)z_j + O(y)\qquad (j=1,\dots, m)\\
\dot u= F(Z,u,v) + O(y), \qquad \dot v =G(Z,u,v), \\
\dot y=H(Z,u,v,y)y.\end{array}
\end{equation}
The center manifold for this map is given by $y=0$ (as it is invariant and is tangent to the eigenspace of the matrix $A$
corresponding to all the multipliers on the unit circle). Thus the restriction of $T$ on $W^c$ equals to $\sigma$
times the time-$1$ map by the flow of
\begin{equation}\label{pnfff00}
\begin{array}{l}\dot z_j=i \Omega_j(Z,u,v)z_j\qquad (j=1,\dots, m)\\
\dot u= F(Z,u,v), \qquad \dot v =G(Z,u,v).\end{array}
\end{equation}

Recall that we have either $h|_{W^c}: (z,u,v) \to (z^*,u,v)$ or $h|_{W^c}: (z,u,v)=(z^*,-u,v)$.
In the first case, the reversibility requires that functions $F,G$ must be identically zero, and $\Omega_j$ must be real.
In the second case, the reversibility allows $F(0,u,0)$ to be an arbitrary even function of $u$. Thus, generically, $F(0,u,0)=au^2+O(u^4)$ with $a\neq 0$,
so by adding an arbitrary constant term to $F$ we can destroy the fixed point. It follows that the only generic case where a symmetric periodic orbit may have
multipliers equal to $+1$ corresponds to $h$ being identity on the $u$-space. The flow normal form (\ref{pnfff00}) then recasts as
\begin{equation}\label{nf}
\dot z_j=i \Omega_j(Z,w)z_j \qquad (j=1,\dots, s),\qquad
\dot w= 0,
\end{equation}
where we denote $w=(u,v)$, and $\Omega_j$ is a real function.
In the polar coordinates $z_j=\sqrt{Z_j}e^{i\phi_j}$
the system takes the form
\begin{equation}\label{nf0}
\dot \varphi=\varphi + \Omega(Z,w), \;\;
\dot Z=0,\qquad \dot w= 0,
\end{equation}
where $\varphi=(\varphi_1,\dots,\varphi_s)$, $\Omega=(\Omega_1,\dots,\Omega_s)$. In these coordinates the involution $h$ is given by $h(\varphi,Z,w)=(-\varphi,Z,w)$.

We further assume that $det\left(\frac{\partial \Omega}{\partial Z}\right)\neq 0$. Then, arbitrarily close to $Z=0$ there exists a value of $Z=Z^*$ such that
the corresponding vector of frequencies $\Omega$ is a rational multiple of $\pi$, i.e. $\Omega(Z^*,0)=\frac{\pi}{q}(p_1,\dots,p_s)$ where $p_j$ and $q>0$ are integers.
It follows that every point in the torus $(Z=Z_j,w=0)$ is $2q$-periodic. Take a non-symmetric periodic orbit on this torus and let $M=(\varphi^*,Z^*,0)$
be a point of this orbit. We can choose the frequencies $\Omega-\Omega(Z^*,0)$ as new $Z$-coordinates near $M$ (so $M$ will be given by $(\varphi=\varphi^*,Z=0,w=0)$
in the new coordinates). The map $T|_{W^c}^{2q}$ in these coordinates is given by
$$\bar \varphi= \varphi+ 2qZ, \qquad \bar Z= Z, \qquad \bar w =\sigma^{2q} w=w.$$
By a small perturbation localized in a small neighborhood of the point $f^{-1} M$ we can bring the map $T|_{W^c}^{2q}$
near $M$ to the form
\begin{equation}\label{mapf}
\bar \varphi= \varphi+ 2qZ, \qquad \bar Z= Z(1-2q\varepsilon)-\varepsilon (\varphi-\varphi^*), \qquad \bar w = w,
\end{equation}
where $\varepsilon>0$ is a small parameter. The perturbation is localized in a small neighborhood of a non-symmetric point $f^{-1}M$,
so by adding an appropriate perturbation localized in a small neighborhood of the point $gf^{-1}M$ we can keep our map in the class of reversible systems;
the map $T|_{W^c}^{2q}$ near the point $M$ would keep its form, as the orbit of $M$ is not symmetric, hence does not enter the
small neighborhood of $gf^{-1}M$. Now note that map (\ref{mapf}) is, at each level $w=const$, a linear rotation near $M$. Therefore, for an appropriate choice of $\varepsilon$, all the points
in a neighborhood of $M$ are periodic with the same period, i.e., we have created a periodic spot on $W^c$ in the small neighborhood of the periodic point $x_0$.

It is shown in \cite{T15} that arbitrarily close in $C^r$ to a map
$f$ with a periodic spot there exists a universal map $\hat f$, such that the perturbation $(\hat f-f)$ is supported in the periodic
spot and, moreover, any $C^r$-generic map close enough to $\hat f$ is universal too. In our situation, we have
a periodic spot, which we denote by $Q$, on the center manifold $W^c$; note that it is away from the set $\mbox{Fix}\;(g)$ of the fixed points of $g$. By \cite{T15}, we can perturb the map $f|_{W^c}$ to make it $d$-universal; the perturbation can be
extended outside of $W^c$ in such a way that it will be supported in a small neighborhood of $Q$. We then make the perturbed map $\hat f$ reversible (i.e., we ensure that $\hat f^{-1}=g\circ \hat f \circ g$) by adding an appropriate small perturbation to $f$ in a neighborhood of $gQ$ (we can do it as $Q$ and $gQ$ do not intersect).

Thus, let $\cal U$ be an open set in the space of $g$-reversible $C^r$-diffeomorphisms such that the generic symmetric periodic orbit $P$ under consideration persists for every map $f$ from $\cal U$. We have shown that for any $\delta>0$, for a residual subset ${\cal U}_\delta$ of $\cal U$ the restriction of any map from ${\cal U}_\delta$ to the local center manifold $W^c(P)$ in the $\delta$-neighborhood of $P$ is $d$-universal, where $d=dim(W^c(P))$. Every map which belongs to the intersection of the sets ${\cal U}_{\delta_j}$ over a tending to zero sequence $\delta_j$ is $C^r$-generic in $\cal U$ and
has the desired property: its restriction to $W^c(P)$ in an arbitrarily small neighborhood of $P$ is $d$-universal. $\square$

\section*{Acknowledgments} Results from Sects. 1,2 and 4 were obtained at financial support of the RSF grant 14-41-00044.
Results from Sec. 3 were obtained at financial support of the RSF grant 14-12-00811.
S. Gonchenko thanks Russian Foundation for Basic Research, grant 16-01-00364, and the Russian Ministry of Science and Education, project 1.3287.2017 -- target part, for supporting scientific researches. D. Turaev also acknowledges the support by the Royal Society, EPSRC, and the Imperial College Department of Mathematics Platform Grant.

\end{document}